\documentclass[12pt]{article}
\usepackage{amsmath,amssymb,epic}
\evensidemargin\oddsidemargin
\advance\textheight-\headheight
\advance\textheight-\headsep
\advance\textheight-\topskip

\newtheorem{theorem}{Theorem}

\newtheorem{prop}[theorem]{Proposition}

\newtheorem{defn}[theorem]{Definition}

\newcommand{\proof}{\noindent\underline{Proof.}\ }
\newtheorem{rem}[theorem]{Remark}
\newtheorem{lemma}[theorem]{Lemma}

\renewcommand{\setminus}{-}

\newcommand{\be}{\begin{enumerate}}
\newcommand{\ee}{\end{enumerate}}
\newcommand{\bi}{\begin{itemize}}
\newcommand{\ei}{\end{itemize}}

\newcommand{\ba}{\begin{array}}
\newcommand{\ea}{\end{array}}

\newcommand{\ovl}{\overline}

\newcommand{\lsemi}{\rtimes}
\newcommand{\rsemi}{\ltimes}

\newcommand{\Lie}{\mbox{\rm Lie}}
\newcommand{\SL}{\mbox{\rm SL}}
\newcommand{\GL}{\mbox{\rm GL}}

\newcommand{\red}{\mbox{\scriptsize \rm red}}
\newcommand{\Spin}{\mbox{\rm Spin}}
\newcommand{\Stab}{\mbox{\rm Stab}}
\newcommand{\Hom}{\mbox{\rm Hom}}
\newcommand{\Tr}{\mbox{\rm Tr}}
\newcommand{\ad}{\mbox{\rm ad}}

\newcommand{\zentr}{\mbox{\large{\sf Y}\normalsize}}

\newcommand{\G}{{G}}
\newcommand{\GG}{{\mathbb G}}
\newcommand{\Z}{{\mathbb{Z}}}
\newcommand{\Q}{{\mathbb{Q}}}
\newcommand{\F}{{\mathbb{F}}}

\newcommand{\R}{{\mathbb{R}}}

\newcommand{\C}{{\mathbb{C}}}

\newcommand{\eb}{\phantom{zzz}\hfill{$\square $}\smallskip}

\renewcommand{\em}{\sf}

\newcommand{\knubbel}{\setlength{\unitlength}{0.1cm}\begin{picture}(4,2)(0,0)  \put(1,1){\circle*{2}} \end{picture}}

\title{Globally maximal arithmetic groups.}

\author {
Benedict Gross \\ 
{\small Dept. of Mathematics, Harvard University, One Oxford Street, Cambridge, MA 01238, USA}
\\
{\small gross@math.harvard.edu}
\\
Gabriele Nebe \footnote{Most of this work was carried out during G. Nebe's visit to Harvard University in Spring and Summer of 2002} \\
{\small Abteilung Reine Mathematik, Universit\"at Ulm, 89069 Ulm, Germany} \\
{\small nebe@mathematik.uni-ulm.de} 
}

\begin{document}

\maketitle

\section{Introduction.}

Let $\G $ be a linear algebraic group defined over $\Q $, and assume that 
$\G (\R )$ is compact and meets every connected component of
$\G (\C )$.
Let $\hat{\Q } := \hat{\Z} \otimes \Q $ be the ring of finite ad\`eles.
Every arithmetic subgroup $\Gamma $ of $\G (\Q )$ is finite, and is obtained by
choosing an open, compact subgroup $K$ of $\G (\hat{\Q })$ and defining 
$\Gamma = K \cap \G (\Q )$ in $\G (\hat{\Q })$.
We note that $\G (\Q )$ is discrete and co-compact in $\G (\hat{\Q })$.

In this paper, we consider the cases where the arithmetic subgroup $\Gamma $
is contained in a unique maximal compact subgroup $K_p$
of $\G (\Q _p)$, for all primes $p$.
We call such $\Gamma $ {\em globally maximal};
examples are provided by finite groups $\Gamma $ with globally irreducible
representations $V$ over $\Q $ where $\G $ is a classical group
$O(V)$, $SU_K(V)$, or $SU_D(V)$, according to whether the
commuting algebra of $V$ is $\Q $,  an imaginary quadratic field $K$ or
a definite quaternion algebra $D$.
Other, in general not globally irreducible, examples are provided by the
finite absolutely irreducible 
 rational matrix groups that are
``lattice sparse'' of even type (see \cite{NeP}). These are finite subgroups 
$\Gamma \leq \GL_n(\Q )$
for which the natural representation is absolutely irreducible such that all $\Gamma $-invariant
lattices can be obtained from any $\Gamma $-invariant lattice $L$, by successively taking 
the dual lattice, scalar multiples, intersections and sums of lattices that are already 
constructed (there are many such groups, e.g. for $n=24$ there are 34 such maximal finite
groups).
Here the algebraic group $\G $ is
$\G = O(V)$ 
and the maximal compact subgroup $\G (\Q _p)$ containing
$\Gamma $ is $O(L\otimes \Z_p)$ for any $\Gamma $-invariant lattice $L$.

Another simple example of a globally maximal $\Gamma $ is the group 
$\Gamma = S_4 = 2^2\lsemi \SL_2(2)$,
which has a unique irreducible representation $V$ of dimension 3 and determinant 1.
This representation is orthogonal, and $\Gamma $ is an arithmetic subgroup of $G=SO(V)$.
The unique maximal compact $K_p$ containing $\Gamma $ is hyperspecial, for $p\neq 2$,
and $K_2 = G(\Q _2)$.
In this paper, we will consider similar examples, when $\G $ is the
unique anisotropic form of $G_2$, $F_4$, and $E_8$ over $\Q $.
In these cases, $\G $ is split over $\Q _p$ for all primes $p$.

In \cite{Gross} the first author has already given some examples of
 globally maximal $\Gamma $, where $K_p$ is hyperspecial for all primes $p$.
These are groups over $\Z $, such as $\Gamma = G_2(2)$ in $\G$ of type $G_2$, and
$\Gamma = \ ^3 D_4(2).3 $ in $\G$ of type $F_4$.
Here we will consider the more exotic cases of Jordan subgroups $\Gamma $,
where $K_p$ is not hyperspecial at a single prime $p$.
To identify the maximal parahoric subgroup $K_p$ containing $\Gamma $ at this prime,
we will determine the discriminant of its Lie algebra with respect to a multiple of the 
Killing form.

We begin with a review of the structure of simple, simply-connected complex Lie groups $G=\G(\C )$
and their Lie algebras ${\mathfrak g}_{\C }$. We describe the Chevalley lattice ${\mathfrak g}$
and the associated split group $\G$ over $\Z $.
This gives us a hyperspecial maximal compact subgroup $\G(\Z _p)$ in $\G(\Q_p)$ 
and we describe the other maximal parahoric subgroups $K_p$ and their Lie algebras starting
from $\G(\Z _p)$.
We then consider the Killing form on ${\mathfrak g}$ and show that it is divisible by $2h^{\vee }$,
where $h^{\vee }$ is the dual Coxeter number.
The same holds for the Lie algebras of the other maximal parahorics. We compute the
discriminants of the resulting scaled forms.
Finally we consider the Jordan subgroups $\Gamma=$ 
$2^3.\SL_3(2) \leq G_2$, $3^3\lsemi \SL_3(3) \leq F_4$, 
and $2^5.\SL_5(2) \leq 2^5.2^{10}.\SL_5(2) \leq E_8$ and determine
the $\Gamma $-invariant lattices in ${\mathfrak g}_{\Q }$.
 The $\Gamma $-invariant 
Lie brackets on ${\mathfrak g}$ are unique up to scalar multiples, except for 
$\Gamma =   3^3\lsemi \SL_3(3) \leq F _4$, where there are two possible Lie brackets (which are 
interchanged by an outer automorphism). 
We show that these Jordan subgroups are globally maximal and
determine their maximal compact overgroups $K_p \leq \G (\Q _p)$.
The last section treats the Jordan subgroups of the classical groups.

\section{Simple Lie groups.} (cf. \cite{Serre})

Let $G$ be a simple, simply-connected, complex Lie group.
Let $$T\subset B \subset G $$ be a maximal torus contained in a Borel subgroup
of $G$.
Let $X^{\bullet}$ denote the character group of $T$. 
This is a free abelian group, containing the finite set $\Phi $ of roots - 
the non-zero characters of $T$ which occur on ${\mathfrak g} = \Lie (G)$.
Let $\Phi _+ \subset \Phi $ be the positive roots, which occur on 
$\Lie (B)$, and let $$\Delta \subset \Phi _+$$ be the root basis
determined by $B$.
Every root $\beta $ in $\Phi _+$ can be written uniquely as 
$$\beta = \sum _{\alpha \in \Delta } n_{\alpha }(\beta ) \alpha \ \ n_{\alpha } \geq 0 .$$
Since $G$ is simple, there is a highest root $\beta _0$ with the
property that
$$n_{\alpha }(\beta _0) \geq 1 \ \ \mbox{ for all } \alpha \in \Delta $$
$$n_{\alpha }(\beta _0) \geq n_{\alpha }(\beta ) \ \ \mbox{ for all } \beta\in \Phi _+, \alpha \in \Delta $$
The sum $$h:= 1 + \sum _{\alpha \in \Delta } n_{\alpha }(\beta _0 ) $$
is the Coxeter number of $G$.
Let ${\mathfrak t} := \Lie (T)$. Then, as a representation of $T$,
$${\mathfrak g} = {\mathfrak t} + \bigoplus _{\beta \in \Phi } {\mathfrak g} ^{\beta } $$
where each root space ${\mathfrak g}^{\beta }$ has dimension 1.
The space ${\mathfrak t}^{\beta } := [ {\mathfrak g}^{\beta } , {\mathfrak g} ^{-\beta }]$
has dimension 1 and is contained in ${\mathfrak t}$.
It has a unique basis $H_{\beta }$ which satisfies $\beta(H_\beta )=2$,
here we have identified $\Hom ({\mathfrak t},\C )$ with 
$X^{\bullet} \otimes \C $. 
If $X_{\beta }$ is a basis for ${\mathfrak g}^{\beta }$, there is a unique
basis vector $Y_{\beta }$ for ${\mathfrak g}^{-\beta }$ with 
$$ [ X_{\beta } , Y_{\beta }] = H_{\beta }.$$
Furthermore we have 
$[H_{\beta } , X_{\beta } ] = 2 X_{\beta } $, $[H_{\beta }, Y_{\beta } ] = -2 Y_{\beta } $.
Hence $${\mathfrak g}_{\beta } :=\langle H_{\beta } , X_{\beta }, Y_{\beta}  \rangle $$ 
is a sub-algebra of ${\mathfrak g}$ isomorphic to ${\mathfrak sl}_2$.
 
By Lie's theorem, the homomorphism ${\mathfrak sl} _2 \to {\mathfrak g}$ 
given by the root $\beta $ lifts to a homomorphism of complex Lie groups
$\SL _2  \to G$. The unipotent subgroup $\GG _a \cong 
\left(\begin{array}{cc} 1 & {\star} \\ 0 & 1 \end{array} \right) $ of $\SL _2$
maps to the root group $U_{\beta }$ of $G$, with Lie algebra ${\mathfrak g}_{\beta }$. 
The map of the tori $\GG _m \cong \left(\begin{array}{cc} t & 0 \\ 0 & t^{-1} \end{array} \right) \to T$ is the co-root
$\beta ^{\vee }$ in $X_{\bullet} = \Hom (X^{\bullet} , \Z )$.
Under the identification 
$X_{\bullet} \otimes \C = \Lie (T)$, $\beta ^{\vee }$ maps to the 
vector $H_{\beta }$ in ${\mathfrak t}^{\beta }$.
Since $G$ is simply-connected, the co-roots span $X_{{\bullet}}$, and the simple 
co-roots $\alpha ^{\vee} $, $\alpha \in \Delta $ give a $\Z $-basis.

The Weyl group $W=N_G(T)/T$ acts on $X_{\bullet}$ and $X^{\bullet}$, the pairing
$X_{\bullet} \otimes X^{\bullet} \to \Z $ is $W$-invariant.
Since $G$ is simple, the action of $W$ on $X_{\bullet} \otimes \Q $ is 
(absolutely) irreducible.
Hence there is a unique $W$-invariant pairing 
$$\langle \phantom{c} , \phantom{c} \rangle : X_{\bullet} \otimes X_{\bullet} \to \Z $$
which is even, indivisible, and positive definite.
We have $$\langle \alpha ^{\vee } , \alpha ^{\vee } \rangle = 2 $$
if $\alpha $ is a long root, and 
 $$\langle \alpha ^{\vee } , \alpha ^{\vee } \rangle = 2c $$
if $\alpha $ is a short root (so $\alpha ^{\vee} $ is a long co-root), with $c=2$ or $3$.
The dual Coxeter number $h^{\vee }$ is defined by 
$$h^{\vee } = 1 + \sum _{\alpha \ \mbox{{\scriptsize long}}} n_{\alpha} (\beta _{0})
+\frac{1}{c} \sum _{\alpha \ \mbox{{\scriptsize short}}} n_{\alpha} (\beta _{0}) $$
Here is a table:
$$\begin{array}{|c|c|c|c|} 
\hline
G & h & c & h^{\vee } \\ \hline
A_n & n+1 & 1 & n+1 \\ \hline
B_n & 2n & 2 & 2n-1 \\ \hline
C_n & 2n & 2 & n+1 \\ \hline
D_n & 2n-2 & 1 & 2n-2 \\ \hline
G_2 & 6 &3 &4 \\ \hline
F_4 & 12 & 2 & 9 \\ \hline
E_6 & 12 & 1 & 12 \\ \hline
E_7 & 18 & 1 & 18 \\ \hline
E_8 & 30 & 1 & 30 \\ \hline
\end{array} $$

\section{Integral Theory}\label{INT}
We now modify our notation slightly: the complex Lie groups and Lie algebras
of the previous section will now be denoted by 
$G_{\C }$ and ${\mathfrak g}_{\C }$ to preserve $G$ and ${\mathfrak g}$ for the integral forms.

Chevalley proved that one can choose the basis elements 
$X_{\beta }$ of the root eigenspaces ${\mathfrak g}^{\beta } \subset {\mathfrak g}_{\C }$ so that 
$$[X_{\beta } , X_{-\beta }] = H_{\beta }$$
$$[X_{\beta } , X_{\alpha }] = 0 \mbox{ if $\alpha +\beta \neq 0 $ is not a root}$$
$$[X_{\beta } , X_{\alpha }] = \pm (m+1)
X_{\alpha + \beta } \mbox{ if $\alpha +\beta $ is a root.}$$
Here $m\geq 0$ is the largest integer such that $\beta - m\alpha $ is a root;
an examination of the root systems of rank 2 shows that $m=0,1, $ or $2$.
The abelian subgroup ${\mathfrak g}$ of ${\mathfrak g}_{\C }$ 
spanned by the $H_{\beta }$ and $X_{\beta }$, $\beta \in \Phi $ is 
a Lie order with $\Z $-basis 
$(H_{\alpha }, X_{\beta } \mid \alpha \in \Delta , \beta \in \Phi )$.
This is the Lie algebra of the split, simply-connected group $\G$ over $\Z $ with
complex points $\G (\C ) = G_{\C }$. The group $\G$ is generated by the integral
Torus $T = X_{\bullet} \otimes \GG _m $, and the root subgroups 
$U_{\beta} \cong \GG _a$ with Lie algebras $\Z X_{\beta }$.

The group $\G(\Z_p)$ gives a hyperspecial maximal compact subgroup of $\G(\Q _p)$ for every prime $p$.
This contains the Iwahori subgroup $I_p$, with reduction to the 
Borel $B$ mod $p$. 
We have
$$\Lie (I_p) = \bigoplus _{\alpha \in \Delta } \Z_p H_{\alpha } 
\oplus \bigoplus _{\beta < 0 } \Z_p X_{\beta } 
\oplus \bigoplus _{\beta > 0 } p\Z_p X_{\beta }  .$$
We want to describe the maximal parahoric subgroups of $\G(\Q_p)$
which contain $I_p$. Besides $\G(\Z _p)$, they are indexed 
by the simple roots $\alpha $ in $\Delta $, and the groups 
$$\{ \G(\Z_p ) , \G_{\alpha }(\Z _p) \mid \alpha \in \Delta \} $$
represent the $(l+1)$ distinct conjugacy classes of maximal compact subgroups of $\G(\Q _p)$.

To each simple root
$\alpha \in \Delta $ 
we can associate a maximal parabolic subgroup $P_{\alpha }$
of $\G (\F_p)$, which contains $B$.
Its inverse image $J_{\alpha}$ in $\G(\Z _p)$ has Lie algebra
$$\Lie (J_{\alpha}) = \bigoplus _{\gamma \in \Delta } \Z_p H_{\gamma } 
\oplus \bigoplus _{n_{\alpha}(\beta )\leq 0 } \Z_p X_{\beta } 
\oplus \bigoplus _{n_{\alpha}(\beta )> 0 } p\Z_p X_{\beta }  .$$
$J_{\alpha}$ is a non-maximal parahoric subgroup, and we will see that
$$J_{\alpha} = \G(\Z _p) \cap \G_{\alpha}(\Z_p ).$$
The next theorem follows from Bruhat-Tits theory.

\begin{theorem}{\label{filt}}
Let $\alpha\in \Delta $ be a simple root.
Then there is a maximal compact subgroup $$G_{\alpha} := \G _{\alpha}(\Z_ p) \leq \G (\Q _p)$$ with Lie-algebra
$$ \Lie (G_{\alpha}) = {\mathfrak g}_{\alpha} :=\bigoplus _{\gamma \in \Delta } \Z_p H_{\gamma }
\oplus \bigoplus _{n_{\alpha }(\beta ) = -n } \frac{1}{p} \Z_p X_{\beta }
\oplus \bigoplus _{-n < n_{\alpha }(\beta )\leq 0 } \Z_p X_{\beta }
\oplus \bigoplus _{n_{\alpha }(\beta )> 0 } p\Z_p X_{\beta }  $$
where $n:=n_{\alpha}(\beta _0)$ is the multiplicity of $\alpha$ 
in the highest root $\beta _0$.
The group $\ovl{G_{\alpha }} := \G_{\alpha }(\F _p)$
is a semidirect product
$$\overline{G_{\alpha}} = G_{\alpha }^{\red } \rsemi R(\overline{G_{\alpha }})$$
where $R(\overline{G_{\alpha }})$ is the unipotent radical, and $G_{\alpha }^{\red }
$ is semi-simple, with root system
$$\Phi _{\alpha } = \{ \beta \in \Phi \mid n_{\alpha }(\beta ) \equiv 0 \pmod{n} \} .$$
This root system has simple roots $\Delta \setminus \{ \alpha \} \cup \{ -\beta _0 \} $
with respect to the Borel subgroup reducing to $I_p$.
The unipotent radical $R:=R(\overline{G_{\alpha }})$ is filtered as a
$G_{\alpha }^{\red}$-module, with $n-1$ abelian subquotients $U_i$
$$R = R_1 \supset R_2 \supset \ldots \supset R_n = \{0 \} $$
$$U_i = R_i/R_{i+1} \cong \bigoplus _{n_{\alpha }(\beta ) \equiv i \pmod{n} } \F_p X_{\beta } .$$
\end{theorem}

\proof
The $\Z_p$-lattice ${\mathfrak g}_{\alpha}$ is a sub-Lie order of
${\mathfrak g} \otimes \Q _p$. 
This will be the Lie-algebra of $\G _{\alpha}$.
Indeed, we may define $G_{\alpha}$ by adjoining to $J_{\alpha}$ the elements 
$e_{\beta }(1/p)$ in the root groups $U_{\beta } \otimes \Q _p$, where $\beta $ is 
a root with $n_{\alpha}(\beta ) = -n = n_{\alpha} (-\beta _0)$, and 
$e_{\beta }: \GG _a \to U_{\beta }$ is the isomorphism over $\Z _p$.
This gives a compact subgroup with desired Lie algebra, by the Chevalley relations.
The theory of Bruhat and Tits shows that $\G_{\alpha}$ defines a smooth group 
scheme over $\Z _p$, and describes its special fiber. 
The filtration of $R$ is obtained by looking at the orbits of the 
Weyl group of $G_{\alpha}^{\red }$ on $\Phi $.
\eb

\begin{rem}
{\rm
When $\G$ is simply-laced, each $U_i$ is a minuscule, irreducible
representation of $G_{\alpha}^{\red} $.
In general, there are at most two orbits of the Weyl group of $G_{\alpha}^{red} $ on the
weights in $U_i$, corresponding to the roots $\beta $ of different
lengths with $n_{\alpha}(\beta ) \equiv i \pmod{n} $.
In this case, $U_i$ need not be irreducible, if $p=c$.
}
\end{rem}

\begin{rem}
{\rm
The semi-direct product structure of $\overline{G_{\alpha}} $ gives a Lie-ideal $M$ with
$$L= \Lie (G_{\alpha}) \supset M \supset pL $$
The quotient $L/M$ is isomorphic to the Lie algebra of $G_{\alpha}^{\red}$ and
$M/pL$ has order $p^{\dim (R(\overline{G_{\alpha}})) } $.
}
\end{rem}

\section{An example - the maximal parahorics in $E_8$.}{\label{E8}}
We illustrate the theory of the previous section with a discussion of the 9 conjugacy classes of maximal parahoric subgroups of $E_8$ over $\Q_p$.
For each, we determine $G_{\alpha }^{\red }$, as well as the minuscule representations in the filtration of $R(\overline{G_{\alpha }})$.
The representations $U_i$ are explicitely identified using the description
of their roots in Theorem \ref{filt} with the help of the 
system LIE \cite{LIE}.

The distinct conjugacy classes of maximal parahoric subgroups of $E_8(\Q _p)$
correspond bijectively to the nodes of the extended Dynkin diagram:
$$
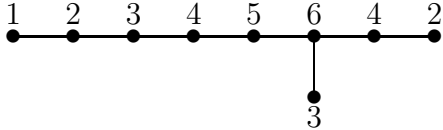
\begin{picture}(200,20)(0,0)
\drawline(10,15)(80,15)
\drawline(60,15)(60,5)
\put(10,15){\circle*{2}}
\put(20,15){\circle*{2}}
\put(30,15){\circle*{2}}
\put(40,15){\circle*{2}}
\put(50,15){\circle*{2}}
\put(60,15){\circle*{2}}
\put(60,5){\circle*{2}}
\put(70,15){\circle*{2}}
\put(80,15){\circle*{2}}
\put(10,17){\makebox(0,0)[cb]{1}}
\put(20,17){\makebox(0,0)[cb]{2}}
\put(30,17){\makebox(0,0)[cb]{3}}
\put(40,17){\makebox(0,0)[cb]{4}}
\put(50,17){\makebox(0,0)[cb]{5}}
\put(60,17){\makebox(0,0)[cb]{6}}
\put(70,17){\makebox(0,0)[cb]{4}}
\put(80,17){\makebox(0,0)[cb]{2}}
\put(60,0){\makebox(0,0)[cb]{3}}
\end{picture}
$$

We have labelled the nodes with the multiplicity 
$n_{\alpha }(\beta _0)$ of the corresponding simple root $\alpha $
in the highest root $\beta _{0 }$.
The extended vertex, with label $n=1$, corresponds to the longest root
$-\beta _0$.

We discuss the parahorics from left to right.
$\mu _a \leq \G_m $ denotes the group of $a$-th roots of unity.
By $\Delta $ we understand a diagonal embedding.
\\
\knubbel
The unique vertex labelled 1 corresponds to the hyperspecial compact
$\G(\Z _p)$. This has 
$$G^{\red} = E_8, R(\overline{G}) = 0, \dim (R(\overline{G})) =0 .$$
\\
\knubbel
The adjacent vertex, labelled 2, has 
$$\begin{array}{ccc}
G_{\alpha }^{\red} & \cong & (\SL_2 \times E_7 )/ \Delta \mu _2 \\
R(\overline{G_{\alpha }}) = U_1 & =  & 2\otimes 56  \\ \dim (R(\overline{G_{\alpha }})) & = &
 112 \end{array}$$
where we have indicated a minuscule representation of a factor by
its dimension. 
\\
\knubbel
The adjacent vertex, labelled 3, has 
$$\begin{array}{ccc}
G_{\alpha }^{\red} & \cong &  (\SL_3 \times E_6 )/ \Delta \mu _3 \\
U_1 & = &  3\otimes 27\\ U_2 & = &  3'\otimes 27' \\ \dim (R(\overline{G_{\alpha }})) & = &  162 \end{array}$$
Here $3'$ is the contragredient representation of the natural representation
$3$ of $\SL_3$, and the representations $27$ and $27'$ of $E_6$ are also dual.
\\
\knubbel
The adjacent vertex, labelled 4, has 
$$\begin{array}{ccc}
G_{\alpha }^{\red} & \cong & (\SL_4 \times \Spin _{10} )/ \Delta \mu _4 \\
U_1 & = & 4\otimes 16 \\ U_2 & = & 6\otimes 10 \\ U_3& = &4'\otimes 16'  \\ \dim (R(\overline{G_{\alpha }})) & = & 188 \end{array} $$
Here $4$ is the natural representation of $\SL_4$, $6 = \Lambda ^2(4)$, and 
$4' = \Lambda ^3(4)$ is the dual of $4$. 
The representations $16$ and $16'$ are the half spin representations of $\Spin(10)$. 
\\
\knubbel
The adjacent vertex, labelled 5, has 
$$\begin{array}{ccc}
G_{\alpha }^{\red}  & \cong  & (\SL_5 \times \SL _{5} )/ \Delta \mu _5 \\
U_1  & = &  5\otimes 10\\ U_2  & = &  10\otimes 5'\\ U_3 & = & 10'\otimes 5 \\
U_4 & = &  5'\otimes 10'\\  \dim (R(\overline{G_{\alpha }}))  & = &  200 \end{array}$$
Here $5$ is the natural representation of $\SL_5$,
$10 = \Lambda ^2 5$, $10' = \Lambda ^3 5$, and $5' = \Lambda ^4 5$.
\\
\knubbel
The adjacent vertex, labelled 6, has 
$$\begin{array}{ccc}
G_{\alpha }^{\red} & \cong  &(\SL_2 \times \SL _{3} \times \SL_6)/ \Delta \mu _6 \\
U_1  & = &  2\otimes 3\otimes 6 \\ U_2  & = &  1\otimes 3'\otimes 15 \\ U_3 & = & 2\otimes 1\otimes 20  \\
U_4 & = &  1\otimes 3\otimes 15' \\ 
U_5 & = &  2\otimes 3'\otimes 6' \\ 
 \dim (R(\overline{G_{\alpha }}))  & = &  202 \end{array} $$
Here $6$ is the natural representation of $\SL_6$,
$15 = \Lambda ^2 6$, $20  = \Lambda ^3 6$, $15' = \Lambda ^4 6$, and
$6' = \Lambda ^5(6) $.
\\
\knubbel
The bottom vertex, adjacent with 6 and labelled 3, has 
$$\begin{array}{ccc}
G_{\alpha }^{\red} & \cong  &\SL_9 / \mu _3 \\
U_1  & = &  84 = \Lambda ^3 9 \\ U_2  & = &  84' = \Lambda ^6 9  \\
 \dim (R(\overline{G_{\alpha }}))  & = &  168 \end{array} $$
where $9$ is the natural representation of $\SL_9$.
\\
\knubbel
The next vertex, adjacent to 6 and  labelled 4, has 
$$\begin{array}{ccc}
G_{\alpha }^{\red} & \cong  &(\SL_2 \times \SL _{8})/ \Delta \mu _2 \\
U_1  & = &  2\otimes 28 \\ U_2  & = &  1\otimes 70 \\ U_3 & = & 2\otimes 28'  \\
 \dim (R(\overline{G_{\alpha }}))  & = &  182 \end{array} $$
Here $8$ is the natural representation of $\SL_8$, $28 = \Lambda ^2 8$, $70  = \Lambda ^4 8$, and $28' = \Lambda ^6 8$.
Similarly $2$ is the natural representation of $\SL_2$ (which is self-dual)
and $1$ is the trivial representation of $\SL _2$.
\\
\knubbel
The last vertex on the right,  labelled 2, has 
$$\begin{array}{ccc}
G_{\alpha }^{\red} & \cong  &\Spin _{16}/ \mu _2 \\
R(\overline{G_{\alpha }})  = U_1  & = &  128 \\ 
 \dim (R(\overline{G_{\alpha }}))  & = &  128 \end{array} $$

In each case it is interesting to note that every 
minuscule representation of $G_{\alpha }^{\red} $ occurs in the
filtration of $R(\overline{G_{\alpha }})$. 
This is a general phenomenon, when $\G$ is of adjoint type, as the center of 
$G_{\alpha }^{\red} $ has order $n= n_{\alpha }(\beta _0)$.

Some other examples of maximal parahorics, which exhibit unusual symmetry, are
given by the following simple roots $\alpha $, indicated in the 
extended Dynkin diagram.

\noindent
$G=\Spin _8:$
$$
\begin{picture}(200,10)(0,0)
\drawline(10,0)(30,10)
\drawline(30,0)(10,10)
\put(10,0){\circle*{2}}
\put(30,0){\circle*{2}}
\put(10,10){\circle*{2}}
\put(30,10){\circle*{2}}
\put(20,5){\circle*{2}}
\put(20,0){\makebox(0,0)[cb]{2}}
\end{picture}
$$
$$\begin{array}{ccc}
G_{\alpha }^{\red} & \cong  &(\SL_2 \times \SL_2 \times \SL_2 \times \SL_2)/ \Delta \mu _2 \\
R(\overline{G_{\alpha }}) = U_1  & = &  2\otimes 2 \otimes 2 \otimes 2\\ 
 \dim (R(\overline{G_{\alpha }}))  & = &  16 \end{array} $$

\noindent
$G=E_6:$
$$
\begin{picture}(200,20)(0,0)
\drawline(10,5)(50,5)
\drawline(30,19)(30,5)
\put(10,5){\circle*{2}}
\put(20,5){\circle*{2}}
\put(30,5){\circle*{2}}
\put(40,5){\circle*{2}}
\put(50,5){\circle*{2}}
\put(30,12){\circle*{2}}
\put(30,19){\circle*{2}}
\put(30,0){\makebox(0,0)[cb]{3}}
\end{picture}
$$
$$\begin{array}{ccc}
G_{\alpha }^{\red} & \cong  &(\SL_3 \times \SL_3 \times \SL_3)/ \Delta \mu _3 \\
U_1  & = &  3\otimes 3 \otimes 3 \\ 
U_2  & = &  3'\otimes 3' \otimes 3' \\ 
 \dim (R(\overline{G_{\alpha }}))  & = &  54 \end{array} $$

\section{The Killing form.}
We retain the notion of Section \ref{INT}, so ${\mathfrak g} = \Lie (G)$ is
the Chevalley Lie algebra of the simply-connected, simple group scheme
$G$ over $\Z$. The Killing form 
$$( X , Y ) := \Tr (\ad X \cdot \ad Y ) $$
is integral, symmetric, and $G$-invariant on ${\mathfrak g}$.
On $X_{\bullet}(T) = \Lie (T)$, it is integral, even, and $W$-invariant,
so it is a multiple of the indivisible form $\langle , \rangle $ with
$\langle \alpha ^{\vee }, \alpha ^{\vee } \rangle = 2$ for $\alpha $ 
a long root.
Steinberg and Springer
\cite{StSp} show that $$(H_{\alpha } , H_{\alpha }) = 4 h^{\vee }$$
for $\alpha $ a long root, with $h^{\vee }$ the dual Coxeter number.
Hence
$$(,) = 2 h^{\vee } \cdot \langle , \rangle $$
as bilinear forms on $\Lie (T)$.
The decomposition $$\Lie (G)  = \Lie (T) \oplus \bigoplus _{\beta > 0 }
(\Z X_{\beta } + \Z X_{-\beta }) $$
is orthogonal for the Killing form. Steinberg and Springer also show that,
for all roots $\beta $, 
$$(X_{\beta }, X_{-\beta }) = \frac{1}{2} (H_{\beta }, H_{\beta }) = h^{\vee} 
\langle H_{\beta } , H_{\beta } \rangle $$
Hence if we define 
$$\langle X, Y \rangle := \frac{1}{2h^{\vee }} (X,Y) $$
we find that
\begin{prop}
The pairing 
$$\langle \phantom{c}, \phantom{c} \rangle : {\mathfrak g} \times {\mathfrak g} 
\to \Z $$ is even, indivisible, and is positive definite on $\Lie(T)$.
\end{prop}
If $G$ is simply-laced, we find that, with respect to 
$\langle \phantom{c}, \phantom{c} \rangle $
$${\mathfrak g}^* / {\mathfrak g} \cong \Lie (T)^* / \Lie (T) \cong \widehat{Z(G) } $$
where $Z(G)$ is the (finite) center of $G$.
In the general case, 
${\mathfrak g}^* / {\mathfrak g} $ has order
$\# {Z(G) } c^k $, where $k$ is the number of short positive roots plus the number of short simple roots.
The latter contribute to $\Lie (T)^* / \Lie (T)$.
\\
Here is a table:
$$\begin{array}{|c|c|c|} 
\hline
G & \det \langle \phantom{c}, \phantom{c} \rangle\mbox{ on } {\mathfrak{g} }  
 & \det \langle \phantom{c}, \phantom{c} \rangle\mbox{ on } {\Lie (T) }   \\
\hline
A_n & n+1 &  n+1 \\ \hline
B_n & 2^{n+2} &  2^2 \\ \hline
C_n & 2^{n^2} &  2^n \\ \hline
D_n & 2^2 & 2^2 \\ \hline
G_2 & 3^7 & 3 \\ \hline
F_4 & 2^{26} & 2^2 \\ \hline
E_6 & 3 &  3 \\ \hline
E_7 &  2 & 2 \\ \hline
E_8 & 1 & 1 \\ \hline
\end{array} $$

The pairing $\langle \phantom{c}, \phantom{c} \rangle$ on 
${\mathfrak g}\otimes \Q _p$ is also integral and even on the $\Z_p$-lattices
$L= \Lie (G_{\alpha })$, for the maximal parahorics in $\G(\Q _p)$ 
defined in Section \ref{INT}.
Indeed the only change in the discriminant $L^*/L$ from that of 
${\mathfrak g}^* / {\mathfrak g} $ involves the planes
$$\Z_p X_{-\beta } + p \Z_p X_{\beta }$$
where $\beta $ is a positive root with 
$$0<n_{\alpha }(\beta ) < n_{\alpha }(\beta _0) .$$
This contributes a factor of $(\Z/p\Z)^2 $ to $L^*/L$. 
Hence we find the following
\begin{prop}
Assume that $p$ does not divide $\det(\langle \phantom{c}, \phantom{c} \rangle )$ on ${\mathfrak g}$.
Then 
$$L^*/L \cong (\Z/p\Z)^{\dim R(\overline{G_{\alpha }})} $$
and $pL^*$ is the Lie ideal $M$ with $L/M \cong \Lie (G_{\alpha }^{\red}) $.
\end{prop}
This allows us to determine which maximal parahorics $\G_{\alpha }(\Z _p)$
can contain certain finite groups $\Gamma \subset \G(\Q _p)$,
once we know some information on the $\Gamma $-stable lattices in 
${\mathfrak g}\otimes \Q_p$.

%
%
%
%
%

\section{The type of some Jordan subgroups of exceptional groups.}
\begin{defn}
Let $\Gamma $ be a globally maximal arithmetic subgroup of a linear algebraic group
$\G $ defined over some number field  $K $.
Then the {\em type} of $\Gamma $ is $${\cal T} (\Gamma ) := ({\cal T}_{\wp }(\Gamma )) _{\wp \mbox{ {\scriptsize prime}}} .$$ Here $\wp $ runs through the
prime ideals of $K$ and
${\cal T}_{\wp }(\Gamma )$ denotes the  maximal compact subgroup of
$\G (K _{\wp })$ over the $\wp $-adic completion of $K$, that contains $\Gamma $.
$$\Gamma \leq {\cal T}_{\wp }(\Gamma ) \leq \G (K _{\wp }).$$
\end{defn}

If $\Gamma $ is a globally maximal group, then the type ${\cal T}_\wp (\Gamma )$ is
hyperspecial for almost all primes $\wp $.
In particular
if the $\Gamma $-module $\Lie (\G (K _\wp ))$
is irreducible modulo $\wp $ then ${\cal T}_{\wp }(\Gamma )  $ is a hyperspecial.

We now treat the different Jordan subgroups of the exceptional groups
in detail.
The explicit calculations are performed using the computer algebra system
MAGMA \cite{MAGMA}.

\subsection{$2^3.\SL_3(2) $ in $G_2$}

The simple roots of $G_2$ are given as follows:
$$
\begin{picture}(100,10)(0,0)
\drawline(10,5)(30,5)
\drawline(20,4.5)(30,4.5)
\put(25,5){\makebox(0,0)[c]{$>$}}
\drawline(20,5.5)(30,5.5)
\put(10,5){\circle*{2}}
\put(20,5){\circle*{2}}
\put(30,5){\circle*{2}}
\put(10,0){\makebox(0,0)[cb]{$-\beta _0$}}
\put(20,0){\makebox(0,0)[cb]{$\alpha _1$}}
\put(30,0){\makebox(0,0)[cb]{$\alpha _2$}}
\end{picture}
$$

The next theorem is already shown in \cite{CNP} by calculations in the 7-dimensional representation of $G_2$.

\begin{theorem}
Let
$\Gamma :=2^3.\SL_3(2)$ be the Jordan subgroup of the anisotropic form 
$\G (\Q )$ of $G_2$.
Then   ${\cal T}_p(\Gamma ) = G_2$ is hyperspecial for $p>2$ and 
${\cal T}_2(\Gamma ) = A_2$.
\end{theorem}

\proof
$\Gamma $ has a unique complex irreducible 14-dimensional representation $V$.
This representation is rational.
The space of $\Gamma $-invariant homomorphisms of $V\otimes V$ to $V$ is one dimensional.
Any generator of this space is skew symmetric and gives a $\Gamma $-invariant Lie-multiplication
on $V$. This yields an embedding of $\Gamma $ into $\G (\Q) $.
The group $\Gamma $ fixes up to isomorphism 12 lattices in $V$ of which the $2$-local inclusions 
are given as follows:

$$\setlength{\unitlength}{0.5mm}
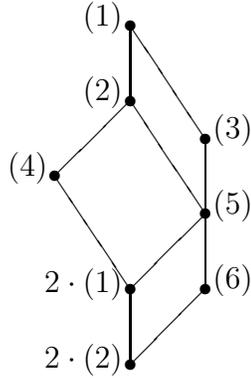
\begin{picture}(100,100)(0,0)
\put(10,60){\circle*{3}}
\put(30,10){\circle*{3}}
\put(30,30){\circle*{3}}
\put(50,30){\circle*{3}}
\put(50,50){\circle*{3}}
\put(50,70){\circle*{3}}
\put(30,80){\circle*{3}}
\put(30,100){\circle*{3}}
\drawline(30,10)(30,30)(10,60)(30,80)(30,100)(50,70)(50,30)(30,10)
\drawline(30,80)(50,50)(30,30)
\put(28,98){\makebox(0,0)[rb]{$\small (1)$}}
\put(28,78){\makebox(0,0)[rb]{$\small (2)$}}
\put(8,58){\makebox(0,0)[rb]{$\small (4)$}}
\put(52,68){\makebox(0,0)[lb]{$\small (3)$}}
\put(52,48){\makebox(0,0)[lb]{$\small (5)$}}
\put(52,28){\makebox(0,0)[lb]{$\small (6)$}}
\put(28,27){\makebox(0,0)[rb]{$2\cdot (1) $}}
\put(28,7){\makebox(0,0)[rb]{$2\cdot (2) $}}
\end{picture} $$

Here the vertical line (e.g from $(1)$ to $(2)$) and the lines parallel to $(2),(4)$ indicate 
inclusions of index $2^3$ (two different $\F_2\Gamma $-modules) and 
the lines parallel to $(1),(3)$ mean inclusions of index $2^8$.
The other 6 isomorphism classes are represented by sublattices of index $3^7$ of these 6 lattices.
This gives the type of $\Gamma $ for all primes $p>2$.
The $\G _{\alpha _1}(\Z_2)$-composition factors of 
 $\Lie (\G_{\alpha _1}(\Z_2))/ 2\Lie( \G_{\alpha _1}(\Z_2))$ are of dimension 1,2, and 4.
Hence the 2-local type of $\Gamma $ is either 
$G_2$ or $A_2$. 
It follows from the mass-formula (see \cite{Gross}, \cite{CNP}) or from the 
calculation in \cite{CNP} that ${\cal T}_2(\Gamma ) = \G _{\alpha _2}(\Z _2)$
is of typs $A_2$.
\eb

\begin{rem}
{\rm
The reduction map  $\Gamma \to \G _{\alpha _2}(\F _2) $ is injective.
}
\end{rem}

\begin{rem}
{\rm
The possibility that ${\cal T}_{2}(\Gamma ) = G_2$ cannot be ruled out looking at the
Lie bracket: 
The maximal $\Gamma $-invariant Lie-order (which corresponds to the lattice (1) in the
picture above) has discriminant $3^7$ (with respect to $1/8$ times the Killing form)
which is the same discriminant as the one of $\Lie (\G (\Z )) $. 
Indeed this Lie-order is also invariant under the maximal compact $\G_{\alpha _2}(\Z _2)$ of type $A_2$ that contains $\Gamma $.
The Lie-order $\Lie (\G_{\alpha _2}(\Z _2))$ corresponds to the lattice number (2) in
the picture above, which is contained in (1) of index $2^3$.
}
\end{rem}

\subsection{$3^3\lsemi \SL_3(3) $ in $  F_4$}
Let $\Gamma $ be the Jordan subgroup $3^3\lsemi \SL_3(3)$ 
 of the unique anisotropic form
$\G (\Q )$ of the algebraic group $F_4$.

The group $\Gamma $ has 3 absolutely irreducible representations of degree 52.
To decide which one is the action of $\Gamma $ on the Lie-algebra of $\G (\Q )$, 
we note that the elements of order 9  in both conjugacy classes of $\G (\Q )$ 
 have trace 1.
This identifies the representation $V = \Lie (\G (\Q)) $ of $\Gamma $ uniquely.
The space $H:=\Hom _{\Gamma }( \Lambda ^2 V , V)$ is 2-dimensional.
The Jacobi identity gives a quadratic equation which has two solutions in $H$. 
Hence there are up to scalar multiples two $\Gamma $-invariant Lie brackets on $V$.
They are interchanged by the outer automorphism 
(in $3^3\lsemi \GL _3(3)$) of $\Gamma $ (which is not in $\G (\Q)$),
 therefore there are up to conjugacy 
two representations of $\Gamma $ into $\G (\Q )$ giving the same conjugacy class of groups 
$\Gamma \leq \G (\Q )$.
We fix  one of the two $\Gamma $-invariant Lie brackets.

The simple roots of $F_4$ are indicated in the following diagram:
$$
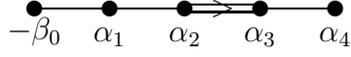
\begin{picture}(100,7)(0,0)
\drawline(10,5)(30,5)
\drawline(30,4.5)(40,4.5)
\put(35,5){\makebox(0,0)[c]{$>$}}
\drawline(30,5.5)(40,5.5)
\drawline(40,5)(50,5)
\put(10,5){\circle*{2}}
\put(20,5){\circle*{2}}
\put(30,5){\circle*{2}}
\put(40,5){\circle*{2}}
\put(50,5){\circle*{2}}
\put(10,0){\makebox(0,0)[cb]{$-\beta _0$}}
\put(20,0){\makebox(0,0)[cb]{$\alpha _1$}}
\put(30,0){\makebox(0,0)[cb]{$\alpha _2$}}
\put(40,0){\makebox(0,0)[cb]{$\alpha _3$}}
\put(50,0){\makebox(0,0)[cb]{$\alpha _4$}}
\end{picture}
$$

\begin{theorem}
${\cal T}_p(\Gamma ) = \G (\Z _p)$ is hyperspecial for $p\neq 3$ and 
${\cal T}_3(\Gamma )  = \G _{\alpha _2}(\Z_3)$.
\end{theorem}

\proof
For $p>3$, the theorem follows because the representation of $\Gamma $ on the
Lie algebra is irreducible modulo $p$.
For $p=2$, this representation 
has two $2$-modular constituents of degree 26.
This implies that ${\cal T}_2(\Gamma )= \G (\Z _2)$ is also hyperspecial.
It remains to consider the prime 3.
The unique maximal $\Gamma $-invariant Lie-order has discriminant $2^{26} \cdot 3^{36}$
(with respect to $1/18$ times the Killing form) which is the discriminant
of the Lie-order $\Lie (\G_{\alpha _2 } (\Z_3))$.
There is no other $\Gamma $-invariant Lie-order that has the discriminant of the Lie algebra
of a maximal compact subgroup of $\G (\Q _3)$.
Therefore
${\cal T}_3(\Gamma ) =\G _{\alpha _2}(\Z_3)$.
\eb

\subsection{$ 2^5.\SL_5(2) $ and $2^5.2^{10}.\SL_5(2)$ in $ E_8$}

Let $\Gamma :=2^5.\SL_5(2) $ be a Jordan subgroup  of the unique anisotropic form
$\G (\Q )$ of the algebraic group $E_8$
and let $H:=2^5.2^{10}.\SL_5(2)$ be the maximal finite Jordan subgroup of $\G (\Q) $ 
that contains $\Gamma $.

The 248-dimensional representation $V$ of $\Gamma $ can be obtained from the
248-dimensional integral representation of the Thompson group, which
contains $\Gamma $ as a maximal subgroup, from the matrices in the www-atlas
\cite{wwwatlas}.
To construct the $\Gamma $-invariant Lie-multiplication on $V$ (which is unique 
up to scalar multiples) we decompose $V$ as the direct sum of eigenspaces 
$$V = \oplus _{\chi } V_{\chi }$$
under the 
maximal normal $2$-subgroup $T\cong 2^5$ of $\Gamma $.
All  31 nontrivial characters $\chi $ of $T$ occur on $V$ with multiplicity 8
and $\Gamma $ permutes the $V_{\chi }$  2-transitively.
The space of $\Stab _{\Gamma }(\chi_1,\chi_2) $-invariant homomorphisms
from $V_{\chi _1} \otimes  V_{\chi _2} $ to $V_{\chi _1\chi_2 }$ is 
onedimensional. From this one constructs the $\Gamma $-invariant Lie-bracket
on $V$.
The maximal decomposable sublattice $L_{OD} := \oplus _{\chi } (V_{\chi } \cap \Lambda )$
of the Thompson-Smith lattice $\Lambda $ carries a $\Gamma $-invariant 
integral Lie-multiplication
such that the discriminant of $L_{OD}$ (with respect to $1/60$ times the Killing form) is
$2^{248}$.

The simple roots of $\G $ are labelled as in the following extended
Dynkin diagram:
$$
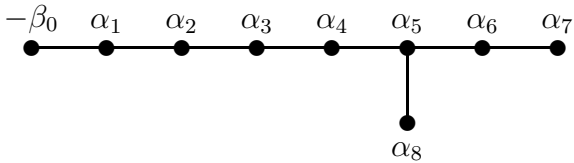
\begin{picture}(200,20)(0,0)
\drawline(10,15)(80,15)
\drawline(60,15)(60,5)
\put(10,15){\circle*{2}}
\put(20,15){\circle*{2}}
\put(30,15){\circle*{2}}
\put(40,15){\circle*{2}}
\put(50,15){\circle*{2}}
\put(60,15){\circle*{2}}
\put(60,5){\circle*{2}}
\put(70,15){\circle*{2}}
\put(80,15){\circle*{2}}
\put(10,17){\makebox(0,0)[cb]{$-\beta_0$}}
\put(20,17){\makebox(0,0)[cb]{$\alpha _1$}}
\put(30,17){\makebox(0,0)[cb]{$\alpha _2$}}
\put(40,17){\makebox(0,0)[cb]{$\alpha _3$}}
\put(50,17){\makebox(0,0)[cb]{$\alpha _4$}}
\put(60,17){\makebox(0,0)[cb]{$\alpha _5$}}
\put(70,17){\makebox(0,0)[cb]{$\alpha _6$}}
\put(80,17){\makebox(0,0)[cb]{$\alpha _7$}}
\put(60,0){\makebox(0,0)[cb]{$\alpha _8$}}
\end{picture}
$$

\begin{theorem}
$\Gamma $ (and hence also $H$) is a globally maximal subgroup of $E_8$.
The type of $\Gamma $  is ${\cal T}(\Gamma ) = {\cal T}(H)$ 
with ${\cal T}_2 (\Gamma )= \G _{\alpha _4}(\Z_ 2)$ and
${\cal T}_p(\Gamma ) = \G (\Z_p)$ for all primes $p>2$.
\end{theorem}

\proof
The 248-dimensional representation of $\Gamma $ is absolutely irreducible modulo 
every prime $p>2$. Therefore ${\cal T}_p(\Gamma ) = {\cal T}_p(H) = \G (\Z_p)$ for
all primes $p>2$.
Since $\Gamma $ is compact, it embeds into at least one maximal compact subgroup $G_{\alpha }$
of $\G (\Q_2)$.
Then $\Gamma $ acts on  $\Lie (G_{\alpha })$
 and hence there is a $\Gamma $-invariant Lie-order
in $V$, of the correct discriminant (with respect to $1/60 $ times the Killing form)
$2^{\dim (R(\overline{G_{\alpha }}))}$ (see Section \ref{E8}).
With MAGMA (\cite{MAGMA}) one calculates that $\Gamma $ fixes up to isomorphism 383 lattices in $V$.
The only lattice of one of the discriminants above, that is closed under the Lie bracket
is a lattice $L_{A4+A4}$ of discriminant $2^{200}$.
Hence $\alpha = \alpha _4$, $\Gamma \leq \G _{\alpha _4}(\Z_2)$, and
${\cal T}_2 (\Gamma ) = \G _{\alpha _4} (\Z _2) $.
\eb

\begin{rem}
{\rm
The representation of $\Gamma $ on $L_{A_4+A_4} / 2 L_{A_4+A_4} $ 
and hence also the reduction map of $\Gamma $ to $\G _{\alpha _4}(\F_2) $  is injective.
From the action of $\Gamma $ on $\Lie (G^{\red} _{\alpha _4})$
one sees that
the image is diagonally embedded in 
$G ^{\red }_{\alpha _4 } \cong \SL_5(2) \times \SL_5(2) / \Delta \mu _5 $.
}
\end{rem}

\begin{rem}
{\rm
There are two maximal $\Gamma $-invariant lattices that are closed under the Lie bracket,
one of which is the orthogonal decomposition $L_{OD} \cong \ ^{(2)}E_8^{31} $ of
discriminant $2^{248}$ and the other lattice is $L_{A4+A4}$.
Therefore $\Gamma $ has 2 maximal Lie orders.
Since both of them are also stable under $H$, the same holds for $H$.
For both Lie orders, the Lie bracket is surjective.
The intersection 
$L_{A4+A4} \cap L_{OD }$ is of index $2^5$ in $L_{OD}$ 
(and of index $2^{29}$ in $L_{A4+A4}$).
}
\end{rem}

\section{The Jordan subgroups of the classical groups.}

In this section we calculate the type of the Jordan subgroups of 
the classical groups $G$. 
To this aim we calculate in the natural representation of $G$.
Then the group $\Gamma $ may contain a center, that acts 
trivially in the adjoint representation of $G$.

\subsection{$p^{1+2n}_+ .  Sp_{2n}(p) \leq U _{p^n}$} 

Let $\Gamma = p^{1+2n}_+ \lsemi  Sp_{2n}(p) \leq U _{p^n}(\C )$ if $p>2$ and 
 $\Gamma  =C_4\zentr 2^{1+2n} . Sp_{2n}(2) \leq U _{2^n}(\C )$ 
if $p=2$.
The minimal number field $L$, such that $\Gamma $ is contained in $U _{p^n}(L)$
is $L=\Q [\zeta_8] $ for $p=2$ and $L = \Q [\zeta_p]$ for $p>2$ and 
the involution is complex conjugation. 
Let $K$ denote the totally real subfield of $L$.
Then the algebraic group $G$ is defined over $K$.

\begin{theorem}
(a) For $p=2$ the type of $\Gamma $ is hyperspecial for all primes
$\wp $ of $K$. 
\\
(b) For $p>2$ the type of $\Gamma $ is hyperspecial for all primes 
 of $K$ not dividing $p$. 
Let $q:=1/2(p^n-1)$. Then 
at the prime $\wp:=(1-\zeta _p)(1-\zeta _p^{-1})$,
the type ${\cal T}_{\wp } (\Gamma )$ is the maximal parahoric
subgroup corresponding to the vertex number $\lfloor \frac{q}{2} \rfloor +1$ 
of the local Dynkin diagram $C-BC_{q}$ (with $q+1$ vertices) (7th diagram on
page 60 \cite{Tits}):
$$
\begin{picture}(100,15)(0,0)
\drawline(20,10)(30,10)
\dashline{3}(30,10)(50,10)
\drawline(50,10)(60,10)
\drawline(60,10.5)(70,10.5)
\drawline(60,9.5)(70,9.5)
\put(65,10){\makebox(0,0)[c]{$>$}}
\drawline(10,10.5)(20,10.5)
\drawline(10,9.5)(20,9.5)
\put(15,10){\makebox(0,0)[c]{$>$}}
\put(10,10){\circle*{2}}
\put(20,10){\circle*{2}}
\put(30,10){\circle*{2}}
\put(50,10){\circle*{2}}
\put(60,10){\circle*{2}}
\put(70,10){\circle*{2}}
\end{picture}
$$
\end{theorem}

\proof
(a) For $p=2$ the group $\langle \zeta _8 \rangle \zentr \Gamma $ is the complex Clifford group described 
in \cite[Section 6]{cliff1}. It follows from the remarks before 
\cite[Theorem 6.5]{cliff1}, that the natural representation of $\Gamma $ is
irreducible modulo all primes $\wp $ of $K$.
Hence the type of $\Gamma $ is hyperspecial everywhere.
\\
(b)
Since $O_p(\Gamma )$ acts absolutely irreducible, the natural representation 
of $\Gamma $ is clearly irreducible modulo all primes of $K$ 
that do not divide $p$. 
\cite{BRW}  shows that the natural representation of $\Gamma $ 
has two $p$-modular constituents of degree $q$ and $q+1$,
where $p^n = 2q+1$. 
Hence $\Gamma $ embeds into the maximal parahoric subgroup $P$ of
 $U _{p^n}(K_{\wp })$ with 
$P^{\red} = O_{q} (p) \times Sp _{q+1} (p)$ if $q$ is odd and 
$P^{\red} = O_{q+1} (p) \times Sp _{q} (p)$ if $q$ is even.
\eb

Note that also for the case $p>2$, the group $\Gamma $ is 
(up to certain scalars) the Clifford group described in 
\cite[Section 7]{cliff1}. 
In particular $\Gamma $ is (up to scalars) a maximal finite 
subgroup of $U_{p^n}(\C )$ (see \cite[Theorem 7.3]{cliff1} for $p>2$
and \cite[Theorem 6.5]{cliff1} for $p=2$).

\subsection{$2^{2n} \lsemi  S_{2n+1} \leq O _{2n+1}$} 

Let $n\geq 3$ and $\Gamma := 2^{2n} \lsemi  S_{2n+1} $. 
Then $\Gamma $ is the determinant 1 subgroup of the full monomial
subgroup $M \leq O_{2n+1}(\Q )$. $M$ is generated by $-I_{2n+1}$
and $\Gamma $. Hence $\Gamma $ fixes the same lattices as
$M$, namely the standard lattice 
$S:=\Z ^{2n+1}$ with quadratic form 
$\sum _{i=1}^{2n+1} x_i^2$, its even sublattice $L$ and the dual lattice $L^*$
(see e.g. \cite{feit}).
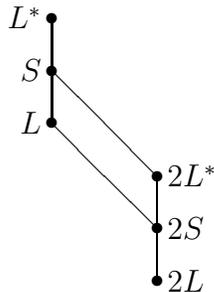
$$\setlength{\unitlength}{0.7mm}
\begin{picture}(100,60)(0,0)
\put(10,60){\circle*{2}}
\put(10,50){\circle*{2}}
\put(10,40){\circle*{2}}
\put(30,30){\circle*{2}}
\put(30,20){\circle*{2}}
\put(30,10){\circle*{2}}
\drawline(10,60)(10,40)
\drawline(10,50)(30,30)
\drawline(10,40)(30,20)
\drawline(30,30)(30,10)
\put(8,58){\makebox(0,0)[rb]{$\tiny L^*$}}
\put(8,48){\makebox(0,0)[rb]{$\tiny S$}}
\put(8,38){\makebox(0,0)[rb]{$\tiny L$}}
\put(32,28){\makebox(0,0)[lb]{$\tiny 2L^*$}}
\put(32,18){\makebox(0,0)[lb]{$\tiny 2S$}}
\put(32,8){\makebox(0,0)[lb]{$\tiny 2L$}}
\end{picture} $$

The unique maximal parahoric subgroup $P$ of $O_{2n+1}(\Q _2, F)$
that contains $\Gamma $ is the orthogonal group of the lattice $L\otimes \Z_2$.
It also stabilizes the dual lattice $L^*$ and, since $L^*/L \cong C_4$,
the unique lattice $S$ between $L^*$ and $L$.
$P$ corresponds to the last vertex of the local Dynkin diagram
$$
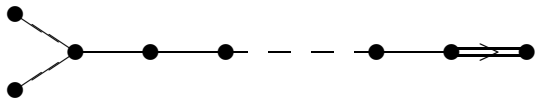
\begin{picture}(100,15)(0,0)
\drawline(2,5)(10,10)(30,10)
\drawline(2,15)(10,10)
\dashline{3}(30,10)(50,10)
\drawline(50,10)(60,10)
\drawline(60,10.5)(70,10.5)
\drawline(60,9.5)(70,9.5)
\put(65,10){\makebox(0,0)[c]{$>$}}
\put(2,5){\circle*{2}}
\put(2,15){\circle*{2}}
\put(10,10){\circle*{2}}
\put(20,10){\circle*{2}}
\put(30,10){\circle*{2}}
\put(50,10){\circle*{2}}
\put(60,10){\circle*{2}}
\put(70,10){\circle*{2}}
\end{picture}
$$

Note that $\overline{P}$ is isomorphic to 
$2^{2n} \lsemi  O_{2n}^+(2)$, if $2n+1 \equiv \pm 3 \pmod{8}$ i.e. the 
2-adic quadratic space is split over $\Q_2$  
(Dynkin diagram $B_n$ on page 60 \cite{Tits}) and 
$\overline{P} \cong 2^{2n}\lsemi O_{2n}^-(2)$, if $2n+1 \equiv \pm 1 \pmod{8}$ i.e. 
the 2-adic quadratic space is non split over $\Q _2$
(Dynkin diagram $\ ^2B_n$ on page 63 \cite{Tits}).

\begin{theorem}
Let $\Gamma := 2^{2n}\lsemi S_{2n+1} \leq O_{2n+1} (\Q )$.
Then $\Gamma $ is globally maximal.
${\cal T}_p(\Gamma )$ is hyperspecial for all primes $p>2$ and
${\cal T}_2(\Gamma ) = P$ as described above.
\end{theorem}

\subsection{$2^{2n-1} \lsemi  S_{2n} \leq O _{2n}$} 
Assume that $n\geq 5$ and let $\Gamma := 2^{2n-1} \lsemi  S_{2n} $
be the determinant 1 subgroup of the full monomial
subgroup $M$ of $O_{2n}(\Q )$.
As in the last section $M$ fixes 3 lattices, 
the standard lattice
$S:=\Z ^{2n}$, its even sublattice $L$ and the dual lattice $L^*$
(see e.g. \cite{feit}).
Since the dimension is even, $L^*/L \cong C_2\times C_2$
and $2L^*$ is contained in $L$.

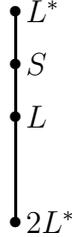
$$\setlength{\unitlength}{0.7mm}
\begin{picture}(100,50)(0,20)
\put(20,70){\circle*{2}}
\put(20,60){\circle*{2}}
\put(20,50){\circle*{2}}
\put(20,30){\circle*{2}}
\drawline(20,70)(20,30)
\put(22,68){\makebox(0,0)[lb]{$\tiny L^*$}}
\put(22,58){\makebox(0,0)[lb]{$\tiny S$}}
\put(22,48){\makebox(0,0)[lb]{$\tiny L$}}
\put(22,28){\makebox(0,0)[lb]{$\tiny 2L^*$}}
\end{picture} $$

Since $\Gamma $ is a normal subgroup of index 2 in $M$, the only other lattices
possibly fixed by $\Gamma $ are the  
 the two lattices $S_1 = \langle L, \frac{1}{2} \sum _{i=1}^{2n} x_i \rangle$ and 
$S_2  =
\langle L, \frac{1}{2} \sum _{i=1}^{2n} x_i - x_1 \rangle $.
These are not fixed under $\Gamma $ (the stabilizer in $M$ of either of these
two lattices is the subgroup of $M$ generated by all permutations and all
even sign changes) and hence $\Gamma $ fixes the same lattices as $M$.

The type of $\Gamma $ is clearly hyperspecial for all primes $p>2$.
For $p=2$, the type $P:={\cal T}_2(\Gamma )$ is as follows:

Assume first that $2n \equiv 0 \pmod{8}$, i.e. the 2-adic quadratic space 
is split. In this case the local Dynkin diagram is 
$D_n$ on page 61 \cite{Tits}.
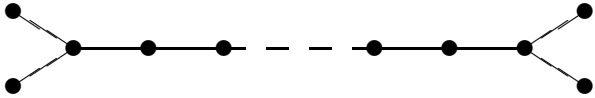
$$
\begin{picture}(100,15)(0,0)
\drawline(78,5)(70,10)
\drawline(78,15)(70,10)
\drawline(2,5)(10,10)(30,10)
\drawline(2,15)(10,10)
\dashline{3}(30,10)(50,10)
\drawline(50,10)(70,10)
\put(2,5){\circle*{2}}
\put(2,15){\circle*{2}}
\put(10,10){\circle*{2}}
\put(20,10){\circle*{2}}
\put(30,10){\circle*{2}}
\put(50,10){\circle*{2}}
\put(60,10){\circle*{2}}
\put(70,10){\circle*{2}}
\put(78,5){\circle*{2}}
\put(78,15){\circle*{2}}
\end{picture}
$$
The two lattices $S_1$ and $S_2$ are even unimodular lattices and their $2$-adic
stabilizers correspond to the two extremal hyperspecial vertices at one side
of the diagram above.
$\Gamma $ interchanges the two lattices $S_1$ and $S_2$ and hence fixes the
midpoint $m$ of the edge joining the two hyperspecial maximal parahorics in the 
building. The type of $\Gamma $ is the maximal compact group $P = Stab (m)$.

Now assume that  $2n \not\equiv 0 \pmod{8}$.
Then the dimension of the anisotropic kernel of the quadratic space is 4
and 
the relative local Dynkin diagram is  $\ ^2 D_n' $ on page 65 \cite{Tits}.

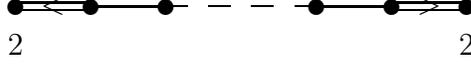
$$
\begin{picture}(100,15)(0,0)
\drawline(20,10)(30,10)
\dashline{3}(30,10)(50,10)
\drawline(50,10)(60,10)
\drawline(60,10.5)(70,10.5)
\drawline(60,9.5)(70,9.5)
\put(65,10){\makebox(0,0)[c]{$>$}}
\drawline(10,10.5)(20,10.5)
\drawline(10,9.5)(20,9.5)
\put(15,10){\makebox(0,0)[c]{$<$}}
\put(10,10){\circle*{2}}
\put(10,5){\makebox(0,0)[c]{$2$}}
\put(70,5){\makebox(0,0)[c]{$2$}}
\put(20,10){\circle*{2}}
\put(30,10){\circle*{2}}
\put(50,10){\circle*{2}}
\put(60,10){\circle*{2}}
\put(70,10){\circle*{2}}
\end{picture}
$$

In this case, $P$ corresponds to one of the two (special) extremal vertices
labelled by 2.

\begin{theorem}
Let $\Gamma := 2^{2n-1} \lsemi  S_{2n} \leq SO _{2n}(\Q )$.
Then $\Gamma $ is  globally maximal.
${\cal T}_p (\Gamma )$ is hyperspecial for $p>2$.
For $p=2$, ${\cal T}_2(\Gamma ) = P$ as described above.
\end{theorem}

\subsection{$2^{1+2n}_+ . O_{2n}^+(2) \leq O _{2^n}$} 

Let $n\geq 3$ and 
$\Gamma := 2^{1+2n}_+ . O_{2n}^+(2) \leq O _{2^n}(\R )$.
Then $\Gamma $ is the full normalizer of the extraspecial group 
$2^{1+2n}_+ $ in the orthogonal group $ O _{2^n}(\R )$.
This gives an isomorphism of $\Gamma $ with the
real Clifford group described in \cite{cliff1}.
Hence up to conjugacy $\Gamma \leq O_{2^n}(\Q [\sqrt{2}]) $ and 
by \cite[Lemma 5.4]{cliff1} $\Gamma $ fixes only one 
lattice in $\Q [\sqrt{2}]^{2^n} $. Hence the type of $\Gamma $ 
is hyperspecial for all primes $\wp $ of $\Z [\sqrt{2}]$.
The reduction modulo $\sqrt{2}$  of $\Gamma $ is the natural 
action of $O_{2n}^+(2) $ on the simple module of the Clifford 
algebra of the associated quadratic form.

\begin{theorem}
The type of 
$\Gamma = 2^{1+2n}_+ . O_{2n}^+(2) \leq O _{2^n}(\Q [\sqrt{2}] )$
is hyperspecial for all primes $\wp $.
\end{theorem}

Note that $\Gamma $ is a maximal finite subgroup of $O_{2^n}(\R )$
as shown in \cite[Theorem 5.6]{cliff1}.

\subsection{$2^{1+2n}_- . O_{2n}^-(2) \leq Sp  _{2^n}$} 

The group 
$\Gamma := 2^{1+2n}_- . O_{2n}^-(2)$ is the centralizer of one factor $Q_8$
in $2^{1+2(n+1)}_+ . O_{2n+2}^+(2) $.
Let ${\cal Q}$ be the quaternion algebra with center $\Q [\sqrt{2}]$ 
ramified only at the two infinite places.
Then $\Gamma $ can be realized as a subgroup of 
$U _{2^{n-1}} ({\cal Q})  \leq Sp _{2^n}(\C )$.

\begin{lemma}
The enveloping order $\Z \Gamma $ of $\Gamma $ in 
${\cal Q}^{2^{n-1}\times 2^{n-1}} $ is a maximal order.
\end{lemma}

\proof
For $n\leq 2$ the lemma can be checked easily by direct computations.
Assume that $n \geq 3$. Then 
$\Gamma $ contains the tensor product 
$\tilde{S}_4\zentr 2^{1+2(n-1)} _+ . O_{2(n-1)}^+(2)$.
Since $n-1 \geq 2$, the group 
$ 2^{1+2(n-1)} _+ . O_{2(n-1)}^+(2) \leq O_{2^{n-1}}(\Q[\sqrt{2}]) $
spans a maximal order $\cong \Z[\sqrt{2}]^{2^{n-1}\times 2^{n-1}} $
by \cite[Lemma 5.4]{cliff1}.
Since $\tilde{S}_4$ spans a maximal order in ${\cal Q}$, the lemma
follows by taking tensor products.
\eb

Hence $\Gamma $ fixes only one class of lattices and hence we get

\begin{theorem}
The type of 
$\Gamma = 2^{1+2n}_- . O_{2n}^-(2) \leq 
U _{2^{n-1}} ({\cal Q})  \leq Sp _{2^n}(\C )$
is hyperspecial for all primes $\wp $.
\end{theorem}


\begin{thebibliography}{ABCD}
\bibitem{BRW}
B. Bolt, T.G. Room, G.E. Wall, {\it On the Clifford collineation,
transform and similarity groups.}
J. Austral. Math. Soc. {\bf 2} (1961) 60-79.
\bibitem{CNP}
A.M. Cohen, G. Nebe, W. Plesken,
{\it Maximal integral forms of the algebraic group $G_2$ defined by finite subgroups},
J. Number Theory 72 (1998) 282-308.
\bibitem{feit}
W. Feit, {\it Integral representations of Weyl groups rationally
equivalent to the reflection representation.}
J. Group Theory {\bf 1} (1998) 213-218.
\bibitem{Gross}
B.H. Gross, {\it Groups over $\Z$},  Inventiones Math. {\bf 124} (1996) 263-279
\bibitem{NeP}
G. Nebe, W. Plesken, {\it Finite rational matrix groups}, AMS-Memoir 
 No. 556, vol. 116 (1995).
\bibitem{cliff1}
G. Nebe, E.M. Rains, N.J.A. Sloane, {\it The invariants of the Clifford groups.}
Designs, Codes and Cryptography {\bf 24} (2001) 99-122.
\bibitem{Serre}
J.-P. Serre, {\it
Alg\`ebres de Lie semisimples complexes.} Benjamin, NY, 1966. 
English translation: {\it Complex semisimple Lie algebras.} Springer (2001) 
\bibitem{StSp} 
T. A. Springer, R. Steinberg 
{\it Conjugacy classes.} 
Seminar on Algebraic Groups and Related Finite Groups (The Institute for Advanced Study, Princeton, N.J., 1968/69) pp. 167--266
Lecture Notes in Mathematics, Vol. 131
Springer (1970)
\bibitem{Tits}
J. Tits, {\it Reductive groups over local fields.} 
Proc. Symp. Pure Maths. {\bf 33} (1979) 29-69.
\bibitem{wwwatlas}
An atlas of finite group representations. Rob Wilson et al. 
http://www.mat.bham.ac.uk/atlas/v2.0/
\bibitem{MAGMA}
 The Magma Computational Algebra System
http://magma.maths.usyd.edu.au/magma/
\bibitem{LIE}
 LiE,
A Computer algebra package for Lie group computations
http://wwwmathlabo.univ-poitiers.fr/$\sim $maavl/LiE/
\end{thebibliography}
\end{document}